\documentclass[12pt]{article}%
\usepackage{amsmath,amssymb,amsthm,enumerate,epsfig,graphicx,natbib, bbm, color,amsfonts}
\usepackage{ifthen,latexsym,syntonly}

\usepackage{verbatim}
\usepackage{morefloats}
\usepackage{multibib}
\usepackage{multirow}
\newcites{app}{References} 
\usepackage{geometry}
\usepackage{soul,color}%
\usepackage{cprotect}
\usepackage{paralist}
\usepackage{float}
\usepackage{booktabs}
\usepackage{bbm}
\usepackage{setspace}
\usepackage[unicode]{hyperref}

\allowdisplaybreaks
\newtheorem{thm}{Theorem}

\newtheorem{prop}{Proposition}

\newtheorem{rmk}{Remark}

\newtheorem{assump}{Assumption}
\newenvironment{myassump}[2][]
{\renewcommand\theassump{#2}
	\begin{assump}[#1]}
	{\end{assump}}

 \DeclareMathOperator{\diag}{diag}
\DeclareSymbolFont{largesymbol}{OMX}{yhex}{m}{n}
\DeclareMathAccent{\Widehat}{\mathord}{largesymbol}{"62}

\def\fX {{\mathbf X}}

\def\fX {{\mathbf X}}

\def\fV {{\mathbf V}}
\def\fW {{\mathbf W}}

\def\fF {{\mathbf F}}
\def\ffb {{\mathbf b}}
\def\fR {{\mathbf R}}
\def\fA {{\mathbf A}}
\def\fB {{\mathbf B}}

\def\fI {{\mathbf I}}

\def\fF {{\mathbf F}}
\def\f1{{\mathbf 1}}

\def\fH{\mathbf H}

\def\fR{\mathbf R}

\def\fSigma{\boldsymbol{\Sigma}}

\def\fme{\boldsymbol{\mathcal{E}}}

\def\me{\mathcal E}

\def\fgamma{{\boldsymbol{\gamma}}}

\def\fPsi{{\boldsymbol{\Gamma}}}
\def\fPsi{{\boldsymbol{\gamma}}}

\def\fmu{{\boldsymbol{\mu}}}

\def\fPsi{\boldsymbol{\Psi}}
\def\fb1{\mathbf{1}_n }
\def\qfb1{\mathcal{Q}( \fb1 )}

\newcommand{\ra}[1]{\renewcommand{\arraystretch}{#1}}

\makeatletter
\renewcommand{\tablename}{{\bf Table}}
\renewcommand{\fnum@table}[1]{\normalfont\textbf{\tablename~\thetable}\\}

\renewcommand{\theenumi}{\Alph{enumi}}

 \makeatletter
 \renewcommand{\p@enumii}{\theenumi.}
 \makeatother

\makeatother
\usepackage{indentfirst} 

\makeatletter
\renewcommand{\figurename}{{\bf Fig.}}
\renewcommand{\fnum@figure}[1]{\normalfont\textbf{\figurename~\thefigure.}}
\makeatother

\makeatletter
\newif\ifrhsapp
\rhsappfalse

\newcommand{\@seccntformat@section}[1]{%
  \ifrhsapp
  Appendix
  \else
  \fi
  \csname the#1\endcsname.\quad
}

\newcommand*{\@seccntformat@subsection}[1]{%
  \csname the#1\endcsname.\quad
}

\newcommand*{\@seccntformat@subsubsection}[1]{%
  \csname the#1\endcsname.\quad
}

\let\@@seccntformat\@seccntformat
\renewcommand*{\@seccntformat}[1]{%
  \expandafter\ifx\csname @seccntformat@#1\endcsname\relax
    \expandafter\@@seccntformat
  \else
    \expandafter
      \csname @seccntformat@#1\expandafter\endcsname
  \fi
    {#1}%
}

\renewcommand{\subsection}{\@startsection
  {subsection}{2}{0mm}{-3.25ex \@plus -1ex \@minus -.2ex}{1.5ex \@plus.2ex}{\normalfont\large\itshape}}

\renewcommand{\subsubsection}{\@startsection
  {subsubsection}{2}{0mm}{-3.25ex \@plus -1ex \@minus -.2ex}{1.5ex \@plus.2ex}{\normalfont\itshape}}

\def\appendix{\par
  \setcounter{section}{0}%
  \setcounter{subsection}{0}%
  \rhsapptrue
  \renewcommand\thesection{\Alph{section}}%
}
\makeatother

\begin{document}

\title{Multiple Testing under High-dimensional  Dynamic Factor Models}
\author{Xinxin Yang$^{a}$ and Lilun Du$^{b,}$\footnote{Correspondence to: College of Business, City University of Hong Kong, 83 Tat Chee Avenue, Kowloon Tong, Hong Kong; Email: lilundu@cityu.edu.hk }  \\
$^{a}$INNO Asset Management  and  $^{b}$City University of Hong Kong
}

\maketitle

\begin{abstract}
Large-scale multiple testing under static factor models is widely used to detect sparse signals in high-dimensional data.
However, static factor models are arguably too stringent because they ignore serial correlation, which seriously distorts error rate control in large-scale inference. In this manuscript, we propose a new multiple testing procedure under dynamic factor models that is robust to nonlinear serial dependence. The idea is to integrate a new sample-splitting strategy based on chronological order and a two-pass Fama--Macbeth regression to form a series of test statistics with marginal symmetry properties and then to use these properties to obtain a data-driven threshold. We show that our procedure can control the false discovery rate asymptotically under high-dimensional dynamic factor models. {As a byproduct of independent interest, we establish a new exponential-type deviation inequality 
for the sum of random variables over various functionals of linear and nonlinear processes.} Our numerical results, including a case study on hedge fund selection, demonstrate the advantage of our proposed method over several state-of-the-art methods.
\end{abstract}

\noindent%


\newpage
\section{Introduction}
Due to advances in information technology, high-dimensional data have become increasingly prevalent in fields such as economics, finance, and healthcare. For instance, in environmental monitoring, key indices such as PM2.5, temperature, and humidity are monitored at numerous stations over time \citep{liang2015assessing}. In fund selection, high-dimensional daily returns are often observed over extended periods \citep{barras2010false}. Similarly, in mobile health, daily ecological momentary assessments (e.g., hours of sleep, daily step counts, and mood scores) are collected via mobile apps, while temporal behavioral patterns are tracked using wristbands preassigned to medical interns \citep{kalmbach2018effects}.  Multiple testing, particularly the false discovery rate (FDR) methodology~\citep{benjamini1995controlling}, has been successfully applied to detect patterns/signals in these scenarios \citep{lan2019factor, wang2021dynamic}.


To capture the dependence of high-dimensional data, factor models are commonly used in finance, economics, and statistics~\citep{fan2021recent}. We begin with a description of traditional static factor models. Assume that the $p \times 1$ vector $\fX_t$ follows a linear factor model: 
\begin{equation}
	\fX_t = \fmu+\fB \fF_{t}+\fme_t ,\quad t=1, \ldots, n, \label{eq:YXZ-1}
\end{equation}
where $\fmu=(\mu_1, \cdots, \mu_p)^\top$ is the mean vector,  $\fF_{t}\in\mathbb{R}^{r \times1}$ denotes the vector of factors whose factor loadings are $\fB_o \in \mathbb{R}^{p \times r} $, and $\fme_t$ is the idiosyncratic component. 
{In static factor models, observations are assumed to be serially independent.}  Statistical inference on the mean vector in static factor models~\eqref{eq:YXZ-1} has received considerable attention over the past two decades. Specifically, the focus is on conducting the following $p$ hypothesis tests simultaneously:
\begin{align}\label{hypo}
	H_{0i}: \mu_i=0\quad  v.s. \quad H_{1i}: \mu_i\neq 0,\quad i=1, \ldots, p,
\end{align}
Let $\mathcal{I}_1 = \{i : \mu_i \ne 0\}$ be the alternative set and $\mathcal{I}_0 = \mathcal{I}_1^c$ represent its complement. We assume that the alternative set is sparse. Our goal is to identify $\mathcal{I}_1$ within the framework of the multiple testing problem \eqref{hypo}.


Various testing procedures \citep{leek2008general, friguet2009factor, fan2012estimating, desai2012cross, fan2017estimation} have been proposed to address the hypotheses \eqref{hypo} under the static factor model \eqref{eq:YXZ-1}.  A common approach involves constructing factor-adjusted test statistics and then applying the standard Benjamini and Hochberg's (BH) procedure \citep{benjamini1995controlling} to identify the alternative set.
When factors are observed, \cite{fan2012estimating} introduce a principal factor approximation method to estimate the false discovery proportion (FDP). As an extension, \cite{fan2017estimation} study the impact of latent factors on the testing procedure and establish a general framework such that the FDP can be well approximated. 
\cite{lan2019factor} propose to eliminate the effect of latent factors through principal component analysis (PCA) and to develop a factor-adjusted procedure to control the FDR.  These studies rely on the assumption of normality of the data. For heavy-tailed data, 
\cite{zhou2018new} construct test statistics for each individual hypothesis using a Huber loss-based approach, and propose their dependence-adjusted procedures. 
\cite{chen2020robust} integrate the adaptive Huber regression with the multiplier bootstrap procedure, which achieves good numerical performance for small or moderate sample size. In this line of research,  a factor-adjusted robust multiple testing procedure is proposed by \cite{fan2019farmtest} in the context of existing latent factors. However, the minimization step to solve the Huber loss and the implementation of multiplier bootstrap or cross-validation techniques in the aforementioned works are computationally expensive. 
\cite{giglio2019thousands} propose a modified BH procedure to detect positive $\mu_i$ values, which is shown to be robust to omitted factors and missing data. \cite{feng2023one} modify the symmetrized data aggregation technique proposed by \cite{Du_2021} to make their method robust to weak/strong factors and various error distributions. All of these methods assume that the observations are serially independent.

\subsection{Dynamic Factor Models}
Serial correlation is a major characteristic of economic and financial data and has been widely studied. For example, \cite{fama1986common} 
study the properties of serial correlation in asset returns and idiosyncratic returns.  \cite{demiguel2014stock} show that exploiting the serial correlation of stock returns improves portfolio performance. Empirical studies on serial correlation in asset returns can also be found in \cite{eun1992forecasting}, \cite{campbell1993trading}, \cite{getmansky2004econometric}, and \cite{fu2009idiosyncratic}.

The causal process provides a principled approach to investigating serial correlation through an explicit measure of dependence.  A sequence $\{X_t\}_{t\in\mathbb{Z}}$ is causal if 
\begin{align}\label{causal}
	X_{t} = g_{t}(\varepsilon_t, \varepsilon_{t-1}, \ldots), \quad \{\varepsilon_t\}_{t\in\mathbb{Z}} \text{~ is an \mbox{i.i.d.} sequence},
\end{align}
where $g_{t}$ indicates measurable functions.
The causal process covers a large class of time series models commonly used in statistics, econometrics, and finance. Two prominent examples are the generalized autoregressive conditional heteroskedasticity (GARCH) sequences and Volterra processes. Examples also include functions of linear processes, $m$-dependent processes, and iterated random functions. See Section~\ref{Sec_3} and \cite{Jirak16} for more details. 

We extend the static factor models in \eqref{eq:YXZ-1} by assuming that each element of the factors or idiosyncratic errors is a causal process while accounting for some observed and confounding factors. Specifically, we consider the following dynamic factor model with confounders:
\begin{align}
	\fX_t=&\fmu+\fB_o \fF_{o, t}+\fB_c \fF_{c, t}+\fme_t,\quad t=1, \ldots, n, \label{eq:YXZ}
\end{align}
where  $\fF_{o, t}\in\mathbb{R}^{r_o\times1}$ denotes the vector of observed factors and $\fF_{c, t}\in\mathbb{R}^{r_c\times1}$ is the vector of confounding factors whose factor loadings are $\fB_o \in \mathbb{R}^{p \times r_o} $ and $\fB_c \in \mathbb{R}^{p \times r_c}$, respectively. Furthermore, each element of vectors $\fF_{o, t}$, $\fF_{c, t}$, and $\fme_t$ is a causal process. In Model \eqref{eq:YXZ}, we do not require linearity between $\fF_{c, t}$ and $\fF_{o, t}$. In fund selection, latent factors represent unobservable systematic risks that remain unaccounted for by commonly used linear pricing models. In microarray data, confounding factors may include environmental conditions and genetic background that influence gene expression levels.


According to \cite{doz2012quasi}\footnote{On page 1015, they mention, ``The model is dynamic since we allow weak serial correlations of the common factor and the idiosyncratic components.''}, static factor models become dynamic when serial correlation exists in both common factors and idiosyncratic components. This point is echoed by \cite{breitung2011gls} where factors are assumed to
be weakly correlated and idiosyncratic errors are autoregressive moving averages (ARMA). It is also  explained by \cite{stock2016dynamic}, according to which the conventional formulation of dynamic factor models can be stacked into a static form as in Models \eqref{eq:YXZ} and \eqref{causal}.
Accommodating both factor structure and serial dependence makes dynamic factor models very effective in analyzing large panels of time series data. 
When studying dynamic factor models, the common technique is to require that factors evolve in specific forms, such as ARMA or GARCH models; see  \cite{stock2016dynamic},\cite{forni2004generalized}, and \cite{breitung2011testing}. In our model, the forms of factors and idiosyncratic components can be unspecified, and nonlinear serial dependence of the data is allowed. In addition, the model can accommodate the nonlinear relationship between observed and confounding factors. These advantages make our model more reliable when analyzing real data.

\subsection{Main Results and Contributions}

New challenges arise when we test large-scale multiple hypotheses~\eqref{hypo} under Models~\eqref{eq:YXZ} and \eqref{causal}. As part of the information about latent factors cannot be explained by a linear function of the observed factors, the cross-sectional correlation between the traditional ordinary least squares (OLS) estimator of $\mu_i$ cannot be eliminated. In addition, 
the widely used PCA method can only generate a biased estimator for factor loadings. 
To obtain an unbiased estimator for each mean $\mu_i$, 
we use the PCA method combined with a two-pass Fama--Macbeth regression \citep{fama1973risk}.
Due to the existence of potentially nonlinear serial correlation, new theoretical techniques are needed to analyze the asymptotic properties of such estimators. 
One of the contributions of our manuscript is that we establish a large deviation bound for the average of random variables with nonspecific (linear or nonlinear) serial correlation (Proposition~\ref{Bernstein_ineq}).
Our large deviation bound is much sharper than that of \cite{liu2010asymptotics} in the sense that it is directly related to the $k$-step serial correlation, with $k$ diverging to infinity. This characteristic also makes our large deviation bound more flexible
when the information from serial correlations is considered. Under mild conditions, we show that factor loadings can be consistently estimated up to a matrix rotation in Proposition~\ref{lem:bhat}, and that the estimator of $\mu_i$ converges uniformly to a normal distribution with a standard deviation that explicitly accounts for serial dependence ({Theorem~\ref{thm:clt}}).

When constructing test statistics, traditional $t$-statistics rely on a consistent estimator of long-run variance in the context of dynamic factor models. Classical kernel-based methods, however, require the selection of a bandwidth, which can significantly influence finite-sample performance (see \cite{andrews1991heteroskedasticity} and \cite{newey1994automatic}). To circumvent the need for bandwidth selection, a self-normalized approach is proposed by \cite{shao2015self, shao2010self}, although with low power.
By constructing new test statistics with a marginal symmetry property,  \cite{zou2020new} propose a computationally efficient and scale-free testing procedure based on sample-splitting. This idea is also adopted to handle more complex data/models in \cite{Du_2021, dai2022false, xing2021controlling, ma2023multiple}. This motivates us to leverage the symmetry of the asymptotic distribution of the proposed estimator without relying on long-run variance estimation. Unfortunately, unlike the case where the data are independent, the serial correlation of the data makes the test statistics for each split used in \cite{zou2020new} highly dependent on each other, leading to a severely distorted FDR. To alleviate this issue, we propose to split the data into two parts in chronological order and take the product of the two estimators computed from each split as the ranking statistic. The next key step is to use the empirical distribution of negative test statistics to approximate the number of false discoveries, so that a data-driven threshold can be obtained by controlling the FDR. 
We show that the FDR of our proposed method can be controlled at the prespecified level asymptotically ({Theorem~\ref{thm:test}}). Our technical derivations go beyond those used in \cite{zou2020new}  due to the existence of nonlinear serial correlation. Another feature of our procedure is that it can accommodate heteroskedasticity and non-sparse observations. 
Numerical results, including a case study on hedge fund selection, are presented to demonstrate the superior performance of our proposed method.

\subsection{Organization of the Paper}
The remainder of this paper is organized as follows. In Section \ref{sec:method},  we introduce our new method for sparse signal detection, followed by the derivation of its theoretical properties in Section~\ref{Sec_3}. Two useful extensions of our proposed method are presented in Section~\ref{sec:extension}. The finite-sample performance of our proposed testing procedure is examined in Section \ref{sec:simu}. In Section~\ref{sec:emp}, we apply our method to study the performance of hedge funds.  Concluding remarks are presented in Section~\ref{sec:conclude}. All technical proofs are reported to the appendix.

We present some notations used below.  Let $\fb1=(1, \ldots, 1)^\top$ and $\mathbf{0}$ denote the zero vector or zero matrix. 
For any matrix $\fA$, define $\mathcal{H}(\fA)=\fA(\fA^\top\fA)^{-1}\fA^\top$ and $\mathcal{Q}(\fA)=\fI-\mathcal{H}(\fA)$ as the projection operator and its complement, respectively. For a given matrix $\fA$,  $\widetilde{\fA} = \mathcal{Q}(\fb1 ) \fA$ is used to demean matrix $\fA$ by column. 
Let  $\lambda_{\min}(\fA)$ denote the smallest eigenvalue of matrix~$\fA$.  For any matrix $\fA$, $\|\fA\|$ represents its spectral norm and $\fA_{i, \cdot}(\fA_{\cdot, i})$ is the $i$-th row (column) of $\fA$.  
For any random variable $X$ and $q\geq1$, let $\|X\|_q=(\mathbb{E}|X|^q)^{1/q}$.  
For any positive constant $x$, $[ x ]$ denotes the largest integer less than $x$.
In the following,  $K$ denotes a generic constant whose value does not depend on $p$ and $n$, and may change from line to line.

\section{Methodology}\label{sec:method}
We use two-pass Fama--Macbeth regressions combined with PCA to estimate each $\mu_i$ in Section~\ref{sec:estalpha}, and then establish a multiple testing procedure based on sample-splitting in Section~\ref{Sec_2.2}. 

\subsection{Estimation of $\fmu$}\label{sec:estalpha}

Assume that $\{\fF_{c, t}, \fF_{o, t}\}$ is a stationary process. By least squares projections \citep{hayashi2011econometrics}, we decompose $\fF_{c, t}$ into two uncorrelated terms: 
\begin{equation}\label{eq:ZW}
\fF_{c, t}:= \fgamma + \fPsi \fF_{o, t}+\fW_t, 
\end{equation}
where $\fgamma  \in  \mathbb{R}^{r_c \times 1}$, $\fPsi \in \mathbb{R}^{r_c \times r_o}$, and $\fW_t \in \mathbb{R}^{r_c \times 1}$
\footnote{Specifically,
$\boldsymbol{\gamma} := \{\mathbb{E}(\fF_{c, t}) - \fPsi \mathbb{E}(\fF_{o, t})\}$, 
$\fPsi^\top := ( \mathbb{E}\{\fF_{o, 1}-\mathbb{E}(\fF_{o, 1})\}\{\fF_{o, 1}-\mathbb{E}(\fF_{o, 1})\}^\top )^{-1}
\mathbb{E}( \{\fF_{o, 1}-\mathbb{E}(\fF_{o, 1})\} \fF_{c, 1}^\top)$, and $\fW_t := ( \fF_{c, t} - \fPsi \fF_{o, t} - \{\mathbb{E}(\fF_{c, t}) - \fPsi \mathbb{E}(\fF_{o, t})\} )$.}.
From \eqref{eq:ZW}, the factor loading $\fPsi$ characterizes the linear extent of confounding between $\fF_{c, t}$ and $\fF_{o, t}$. 
By substituting \eqref{eq:ZW} into \eqref{eq:YXZ}, we obtain
\begin{align}\label{eq:total}
\fX_t= \fmu +(\fB_o+ \fB_c\fPsi  )\fF_{o, t}+\fB_c (\fW_t+\fgamma)+\fme_t.
\end{align}
To identify $\fmu$ in the presence of $\fgamma$, we assume that $\fmu$
is sparse; more details can be found in Assumption~\ref{sparse:alpha} in Section~\ref{Sec_3}.
Let $\fX=(\fX_1, \ldots, \fX_n)\in \mathbb{R}^{p \times n}$, $\fF_o=(\fF_{o, 1}, \ldots, \fF_{o, n})^\top \in \mathbb{R}^{n \times r_o}$, $\fW=(\fW_1, \ldots, \fW_n)^\top \in \mathbb{R}^{n \times r_c}$, and $\fme=(\fme_1, \ldots, \fme_n)\in \mathbb{R}^{p \times n}$.
Then the matrix form of \eqref{eq:total} can be expressed as 
\begin{align}\label{eq:matrix}
\fX= \fmu \fb1^\top +(\fB_o+ \fB_c\fPsi )\fF_{o}^\top +\fB_c(\fW+\fb1  \fgamma^\top)^\top +\fme.
\end{align}
If we directly use OLS regression on $\{\fX, \fF_o\}$, by simple computation, we can estimate the mean vector as
\begin{align}\label{eq:orghat_alpha}
\widehat{\fmu}_{ \mathrm{OLS} }  = \fmu + \fB_c\fgamma + (\fB_c\fW^\top +\fme)\mathcal{Q}(\fF_o)\fb1  (\fb1^\top \mathcal{Q}(\fF_o)\fb1)^{-1}.
\end{align}
By \eqref{eq:orghat_alpha}, the presence of latent factors makes $\widehat{\fmu}_{\mathrm{OLS}}$ not only a biased estimator of $\fmu$ but also results in a strong correlation of its elements due to the common factor $\fW$. This estimator will lead to severely distorted FDRs. 
To address this issue, we sequentially regress the observed and latent factors from $\fX$, which involves the following three steps:

\vspace{1em} 
\textbf{Step I: Time series regression for factor loadings $\fB_o+\fB_c \fPsi$.}  As $\fF_o$ is observed, to remove the term $(\fB_o+\fB_c\fPsi)\fF_{o}^\top$ from $\fX$, we need to estimate $\fB_o+\fB_c\fPsi$. Recall that $\widetilde{\fF}_o = \mathcal{Q}(\fb1)\fF_o$. By directly using the OLS regression on $\{\fX, \fF_o\}$, we obtain the estimator of the coefficient $\fB_o+\fB_c\fPsi$, which can be written as 
\begin{equation}\label{bobc}
(\widehat{\fB_o+\fB_c\fPsi})^\top = (\widetilde{\fF}_o^\top\widetilde{\fF}_o)^{-1}\widetilde{\fF}_o^\top\fX^\top.
\end{equation}
Then, by elementary computation,  we obtain the adjusted return matrix as 
\begin{equation}\label{eq:f_oresid}
\begin{aligned}
\fX-(\widehat{\fB_o+\fB_c\fPsi })\fF_o^\top = &\fX (\fI_n - \widetilde{\fF}_o (\widetilde{\fF}_o^\top \widetilde{\fF}_o)^{-1}\fF_o^\top  )  \\
=& \fmu \fb1^\top  + (\fB_c(\fW+\f1_n\fgamma^\top)^\top +\fme)(\fI_n- \widetilde{\fF}_o(\widetilde{\fF}_o^\top \widetilde{\fF}_o )^{-1} \fF_o^\top), \\
=&  (\fmu+\fB_c\fgamma )\fb1^\top + (\fB_c\fW^\top +\fme)(\fI_n- \widetilde{\fF}_o(\widetilde{\fF}_o^\top \widetilde{\fF}_o )^{-1} \fF_o^\top).
\end{aligned}
\end{equation}

We observe that the right hand of \eqref{eq:f_oresid} has two components: the constant term and the unobserved factor structure.  Both components 
include the term $\fB_c$. To obtain an unbiased estimator of $\fmu$ and eliminate the factor structure, we must first estimate $\fB_c$. 

 \vspace{1em} 
\textbf{Step II: {PCA for factor loadings $\fB_c$.}}  Due to the presence of the constant term in \eqref{eq:f_oresid}, we demean both sides of \eqref{eq:f_oresid} by column before applying PCA. Define
\begin{equation*}
	\fR :=  (\fX- (\widehat{\fB_o+\fB_c\fPsi})\fF_o^\top)\mathcal{Q}( \fb1 )  =  \fX\mathcal{Q}(\widetilde{\fF}_o)\mathcal{Q}( \fb1 ) = (\fB_c\fW^\top+\fme) \mathcal{Q}(\widetilde{\fF}_o)\mathcal{Q}( \fb1 ),
\end{equation*}
which only contains the factor structure. If the number of latent factors $r_c$ is known, using PCA, 
the estimator of $\fB_c$ can be written as 
\begin{align}\label{estbc}
\widehat{\fB}_c=(\widehat{\ffb}_1, \ldots,  \widehat{\ffb}_{r_c}),
\end{align}
where $\widehat{\ffb}_i/\sqrt{p}$ is the eigenvector corresponding to the $i$-th largest eigenvalue $\lambda_i^{\fR\fR^\top }$ of $\fR\fR^\top$, $i=1, \ldots, r_c$. If $r_c$ is unknown, we replace $r_c$ by its estimator $\widehat{r}_c$, which can be obtained by maximizing the eigenvalue ratios as $\widehat{r}_c=\text{argmax}_{1\leq i\leq p-1}\lambda_i^{\fR\fR^\top }/\lambda_{i+1}^{\fR\fR^\top }$.
 
 \vspace{1em} 
\textbf{Step III: Cross-sectional regression for $\fmu$.}  To obtain the estimator of $\fmu$, we use the Fama--MacBeth regression on the left-hand side of \eqref{eq:f_oresid} and $\widehat{\fB}_c$, and get 
\begin{equation}
\begin{aligned}\label{alpha_formula}
\widehat{\fmu} =&  n^{-1} \mathcal{Q}(\widehat{\fB}_c )\fX \big( \fI_n- \widetilde{\fF}_o(\widetilde{\fF}_o^\top \widetilde{\fF}_o)^{-1}\fF_o^\top \big)\fb1.
\end{aligned}
\end{equation}
Thanks to our unique projection matrix (i.e., $\fI_n- \widetilde{\fF}_o( \widetilde{\fF}_o^\top \widetilde{\fF}_o )^{-1}\fF_o^\top$),  the components of $\fX$ associated with the observed factors (i.e., $\fB_o\fF_o^\top$) are removed directly due to the identity
\[ \fF_o^\top(\fI_n- \widetilde{\fF}_o( \widetilde{\fF}_o^\top \widetilde{\fF}_o )^{-1}\fF_o^\top) =\mathbf{0}. \]

\begin{rmk}
{From the proof of Proposition \ref{lem:bhat}, we can see that $\widehat{\fB}_c$ is not a consistent estimator for $\fB_c$. 
In fact, there exists a matrix $\fH_n$ such that $\|\fB_c-  \widehat{\fB}_c \fH_n^{-1} \|_{\infty} =O_p(n^{-1/2})$, where $\|\fA\|_{\infty}$ denotes the maximum absolute row sum of any matrix $\fA$. See Proposition \ref{lem:bhat} for the definition of $\fH_n$. } 
\end{rmk}
\begin{rmk}
{In the Fama--MacBeth regression process, we cannot obtain the unbiased estimator for $\fB_c\fgamma$ and $\fB_c\fW^\top(\fI_n-\widetilde{\fF}_o(\widetilde{\fF}_o^\top \widetilde{\fF}_o)^{-1}\fF_o^\top )\fb1$.  Fortunately, we can obtain the consistent estimator for $\fB_c\fgamma+\fB_c\fW^\top (\fI_n- \widetilde{\fF}_o(\widetilde{\fF}_o^\top \widetilde{\fF}_o)^{-1}\fF_o^\top )\fb1$, which can help us immediately eliminate the bias and factor structure of $\widehat{\fmu}_{\mathrm{OLS}}$ in \eqref{eq:orghat_alpha}.}
\end{rmk}


\subsection{Multiple Testing Procedure via Sample-splitting}\label{Sec_2.2}
The estimator of $\fmu$ provides the input to conduct multiple tests. Typically, the asymptotic null distribution of $\widehat{\mu}_i$ should be established to obtain the $p$-value via normal calibration. However, as idiosyncratic errors and factors are generated from causal processes in dynamic factor models, serial correlations may exist. 
Establishing the asymptotic normality of the estimated $\mu_i$ poses a challenge. Another contribution of our manuscript is to demonstrate that the statistic $\sqrt{n}(\widehat{\mu}_i-\mu_i)/\sigma_i$ converges uniformly to a standard normal distribution for $i$, where $\sigma_i^2$ is the long-run variance under the causal process. Details can be found in Section~\ref{Sec_3}.

However, there are some difficulties in using normal calibration. First, long-run variance is difficult to estimate because it involves selecting an appropriate bandwidth \citep{lobato2001testing}. Second, the accuracy of a multiple testing procedure based on normal calibration depends critically on the skewness of the observation and the diverging rate of $p$ relative to the sample size $n$ \citep{Fan_2007, Delaigle_2011}. The first-order inaccuracy of normal calibration due to skewness can be removed via bootstrapping \citep{liu2014phase}. Accordingly, the number of tests $p$ can diverge faster and a better error rate can be achieved when the population distributions are skewed. However, the computational burden required for  bootstrapping is heavy, especially when the number of hypotheses is large.
To balance this computational issue and robustness, \cite{zou2020new} propose a sample-splitting method to construct a  series of new ranking statistics with symmetry properties and use these properties to approximate the number of false rejections nonparametrically. The proposed procedure is computationally efficient, as it uses only one split of the data. In their simulation studies, this method is shown to be robust to nonnormal distributions with finite moments (e.g., $t$ and chi-square distributions). However, the sample-splitting strategy and FDR control theory in \cite{zou2020new} critically depend on the independence assumption, which is violated in dynamic factor models. 

We propose to split the sample into two halves equally in chronological order and then obtain two estimators of $\fmu$ using the three-step method proposed in Section~\ref{sec:estalpha}. To aggregate the evidence from the two splits, the product of the two estimators serves as the new ranking statistic. The idea is that if the two splits do not substantially overlap, the test statistics computed from these two splits are weakly dependent, so the symmetry property of the ranking statistics still holds asymptotically. We develop the above idea in detail below:
 \begin{enumerate}[(i)]
 \item Sort the data $\{\fX_t, \fF_{o, t}\}$ by time point $t$ in descending order;
 \item Divide the data into two parts: the first part consists of the first $[n/2]$ observations, and the second part consists of the remaining observations;
 \item For each part, estimate $\fmu$ using our proposed algorithm in Section \ref{sec:estalpha}, and denote the two estimates by 
 $\widehat{\fmu}^{(1)}=(\widehat{\mu}_1^{(1)}, \ldots, \widehat{\mu}_p^{(1)})^\top$ and $\widehat{\fmu}^{(2)}=(\widehat{\mu}_1^{(2)}, \ldots, \widehat{\mu}_p^{(2)})^\top$;
 \item Let $T_i^{(1)}=\sqrt{n}\widehat{\mu}_i^{(1)}$ and $T_i^{(2)}=\sqrt{n}\widehat{\mu}_i^{(2)}$. The test statistics for the hypotheses in \eqref{hypo} are defined as
 \begin{align*}
 T_i=T_i^{(1)}T_i^{(2)}, \quad i=1, \ldots, p. 
 \end{align*} 
 \end{enumerate}
  
Under the null hypotheses, we show that the asymptotic distributions of $T_i$ are symmetric with zero mean. 
However, when the $i$-th alternative hypothesis is true, the test statistic $T_i$ is likely to be large, of the order of $O_p(n\mu_i^2 )$. Similar to the nonparametric methods proposed in \cite{barber2015controlling, zou2020new, Du_2021}, we choose a threshold $L>0$ by setting
\begin{equation}\label{LL}
L =\inf\left\{ \xi>0: \frac{1+\sum_{i} I(T_i\leq -\xi) }{\sum_i I(T_i\geq  \xi )\vee1}\leq \beta \right\},
\end{equation}
and reject $H_{0 i}$ if $T_i\geq L$, where $\beta$ is the target FDR level. 
If the set is empty, we simply set $L =+\infty$ and do not make any rejections. Essentially, in \eqref{LL} we use the empirical distribution of negative statistics to approximate the number of false rejections. 

\begin{rmk}
When applying the conventional BH procedure and its ramifications, we need to compute the $p$-value of the test statistic  for each hypothesis in \eqref{hypo}, or other analogous quantities. Thus, it is necessary to have exact information about the asymptotic distribution of the test statistics. However, our proposed procedure only uses the symmetry property of the asymptotic distribution. In particular, we avoid estimating the standard error of the estimated $\mu_i$ for each individual test. Although this may create a fairness issue when comparing $\mu_i$ with different precisions, we show that the proposed test statistics still have the symmetry property and that the FDR can be controlled at the prespecified level. In Section~\ref{Sec_4.1}, an estimator of the individual standard error under a causal process is proposed as input to standardize each $\widehat{\mu}_i$. 
\end{rmk}

\begin{rmk}
The threshold $L$ proposed in \eqref{LL} shares a similar spirit to the knockoff \citep{barber2015controlling}, but they differ in that our procedure does not require any prior information about the distribution of high-dimensional idiosyncratic errors. This is especially important because it is difficult to estimate the distribution when $p$ is very large.
\end{rmk}

\section{Theoretical Results}\label{Sec_3}
In this section, we derive the asymptotic distribution of the estimator $\hat{\mu}_i $ and prove that the FDR of our proposed method can be controlled at the prespecified level. Both tasks depend crucially on sharp deviation inequality for causal processes.

\subsection{A New Concentration Inequality for Causal Processes}
We now derive a new concentration inequality for the average of random variables over a variety of functional processes in Proposition~\ref{Bernstein_ineq}, which may be of independent
interest. Recall that any process $\{X_t\}_{t\in\mathbb{Z}}$ is causal if 
\begin{align*}
X_t=g_t(\varepsilon_t, \varepsilon_{t-1}, \ldots, ),  ~~ t\in\mathbb{Z},
\end{align*}
where $g_t$ represents measurable functions and $\{\varepsilon_t\}_{t\in\mathbb{Z}}$ is a zero mean i.i.d. sequence. 
The causal process encompasses a variety of rich classes. For instance, the GARCH(p, q) sequence, which is given by 
\begin{equation*}
X_t = \varepsilon_t L_t \text{~with~} 
L_t^2 = \mu_L + \alpha_1L_{t-1}^2 +\cdots + \alpha_pL_{t-p}^2+\beta_1X_{t-1}^2 + \cdots + \beta_q X_{t-q}^2,
\end{equation*}
where $\mu_L, \alpha_1, \cdots, \alpha_p, \beta_1, \cdots, \beta_q \in \mathbb{R}$, can be represented as 
\begin{equation*}
	X_t = \sqrt{\mu_L}\varepsilon_t\Bigg(1+\sum_{n=1}^\infty\sum_{1\leq \ell_1\cdots\ell_n\leq r}\prod_{i=1}^n\big(\alpha_{\ell_i}+\beta_{\ell_i}\varepsilon_{t-\ell_1-\cdots -\ell_i}^2\big)\Bigg), \text{~~with $r = \max\{p, q\}$}.
\end{equation*}
Volterra processes, which are fundamentally important for nonlinear processes, can be written according to \cite{berkes2014komlos} as  
\begin{equation*}
	X_t = \sum_{i=1}^\infty \sum_{0\leq j_1<\cdots<j_i} a_t(j_1, \cdots, j_i)\varepsilon_{t-j_1}\cdots\varepsilon_{t-j_i},
\end{equation*}
where $a_t$ is called the $t$-th Volterra kernel.

To characterize the concentration bound,  we introduce a measure of dependence for causal processes. For any stationary causal process $X_{t} = g_t(\varepsilon_{t}, \varepsilon_{t-1}, \ldots)$, define  
$$X_t^{(\ell, \prime)} = g_t(\varepsilon_{t}, \varepsilon_{t-1}, \ldots, \varepsilon_{t-\ell}^\prime,  \varepsilon_{t-\ell-1}, \cdots),
$$
where $\{\varepsilon_{t}^\prime\}_{t\in\mathbb{Z}}$ is an independent copy of $\{\varepsilon_{t}\}_{t \in\mathbb{Z}}$. For simplicity, we denote $X_t^{(t, \prime)}$ by $X_t^{\prime}$.  We define $\Theta_{X, d}(k)=\sum_{t = k}^\infty\|X_t-X_t^{\prime}\|_d$.
It is important to note that $\Theta_{X_, d}(k)$ is closely related to $\mathrm{cov}(X_0, X_k)$ by the following inequality: 
\begin{align*}
	\big|\mathrm{cov}(X_0, X_k)\big|\leq 2\|X_0\|_2\sum_{t = k}^\infty\|X_t-X_t^{\prime}\|_2= 2\|X_0\|_2 \cdot \Theta_{X, 2}(k).
\end{align*}
See \cite{wu2009asymptotic} for the detailed derivation.

\begin{prop}\label{Bernstein_ineq}
	Suppose that $\{\varepsilon_t\}_{t\in\mathbb{Z}}$ is a sequence of \mbox{i.i.d.} random variables, and 
	\begin{align*}
		X_t=g_t(\varepsilon_t, \varepsilon_{t-1}, \ldots, ),  ~~ t\in\mathbb{Z},
	\end{align*}
	is a stationary process with $\mathbb{E}(X_t)=0$ and $\mathbb{E}|X_t|^d<\infty$ for some $d\geq4$. Assume that $\sup_{1\leq t\leq n}|X_t|\leq \sqrt{n}$. Then for any positive sequence $\tau_n/\log n\rightarrow\infty$ and $c_n=o(\log([\tau_n/\log n]))$, we have 
	\begin{align}\label{bernts}
		\mathbb{P}\Big(\Big|\frac{1}{n}\sum_{t=1}^nX_t\Big|\geq \frac{\tau_n}{\sqrt{n}}\Big)=O(n^{-c_n})+\frac{K\cdot\Theta_{X, d}^d([\tau_n])}{\tau_n^d},
	\end{align}
	where $K$ is a constant that depends solely on $d$.
\end{prop}

\begin{rmk}
	From \eqref{bernts}, we find that unlike the case where the random variables $X_t$ are independent, the convergence rate of the stationary process depends on both the sample size and the autocorrelations. 
\end{rmk}
\begin{rmk}
	We compare Lemma 1 in \cite{liu2010asymptotics} with Proposition \ref{Bernstein_ineq}. Assume that $\mathbb{E}|X_1|^d<\infty$, with $d>2$. By Lemma 1 in \cite{liu2010asymptotics},  for any $\tau_n=o(\sqrt{n})$,  we have
	\begin{equation}
		\begin{aligned}\label{bound_wu}
			\mathbb{P}\Big(\Big|\frac{1}{n}\sum_{t=1}^nX_t\Big|\geq \frac{\tau_n}{\sqrt{n}}\Big) \leq&\frac{K\cdot\Theta_{X, d}^d(0) }{\tau_n^d}.
		\end{aligned}
	\end{equation}
	Examining Equations \eqref{bernts} and \eqref{bound_wu} reveals that the convergence rate in \eqref{bernts} is much sharper for the truncated case with $\sup_{1\leq t\leq n}|X_t|\leq \sqrt{n}$.
	For the general case where the assumption $\sup_{1\leq t\leq n}|X_t|\leq \sqrt{n}$ is not satisfied, we have,
	by Markov's inequality and Proposition \ref{Bernstein_ineq}, 
	\begin{equation}\label{new_bound_our}
		\mathbb{P}\Big(\Big|\frac{1}{n}\sum_{t=1}^nX_t\Big|\geq \frac{\tau_n}{\sqrt{n}}\Big) \leq 
		\frac{\mathbb{E}|X_1|^d}{n^{d/2-1}} + \frac{K(d)\cdot\Theta_{X, d}^d (\tau_n)}{\tau_n^d}.
	\end{equation}
	Comparing \eqref{new_bound_our} and \eqref{bound_wu}, we find that both terms in \eqref{new_bound_our} are sharper than those in \eqref{bound_wu}  by appropriately choosing some $\tau_n$ for any $d\geq4$. For instance, if we set $\tau_n = n^{1/4}/\log\log(n)$, then the ratio between the bound in \eqref{bound_wu} and that in \eqref{new_bound_our} tends to infinity. Thus, with \eqref{new_bound_our}, we establish that the dimension $p$ can diverge toward infinity at a higher rate than \eqref{bound_wu}.
\end{rmk}

\begin{rmk}
	To prove Proposition \ref{Bernstein_ineq}, we first construct an $m$-dependent approximation based on $\{X_t\}_{t\in\mathbb{Z}}$, and then decompose it into two sequences, each a sum of independent random variables. This step follows the approach in \citep{Jirak16}. Given the assumption that $\sup_{1\leq t\leq n}|X_t|\leq \sqrt{n}$, we can obtain a bound with exponential decay for both sequences by carefully selecting $m$, which gives the term $O(n^{-c_n})$. The term $\frac{\Theta_{X, d}^d([\tau_n])}{\tau_n^d}$ comes from the error term and explains its dependence on $\tau_n$.
\end{rmk}

\subsection{CLT for $\widehat{\mu}_i$}\label{Sec_3.1}
We now study the asymptotic properties of the estimator $\widehat{\mu}_i$ in Equation~\eqref{alpha_formula}. Examining Equation \eqref{alpha_formula}, we find that the key to analyzing the asymptotic behavior of $\widehat{\mu}_i$ is the term $\widehat{\fB}_c$.  However, due to the nonlinear serial correlation of the data $\{\fX_t, \fF_{o, t}\}$, current technical tools for the linear process \citep{breitung2011gls} cannot be used to analyze the asymptotic properties of $\widehat{\fB}_c$, especially in a high-dimensional setting. Thanks to the new concentration inequality for causal processes in Proposition~\ref{Bernstein_ineq}, we can establish the convergence results of $\widehat{\fB}_c $ under some mild assumptions.

Let $\fF_{o, t}=(F_{o, 1, t}, \ldots, F_{o, r_o, t})^\top$, $\fF_{c, t}=(F_{c, 1, t}, \ldots, F_{c, r_c, t})^\top$, and $\fme_t=(\me_{1, t}, \ldots, \me_{p, t})^\top$. Assume that 
\begin{equation}
	\begin{aligned}\label{causal_1}
		F_{o, i, t}=&F_{o, i, t}(\varepsilon_{f_o, i, t}, \varepsilon_{f_o, i, t-1}, \ldots, ),~~\{\varepsilon_{f_o, i, t}\}_{t\in\mathbb{Z}} \text{~ is an \mbox{i.i.d.} sequence for each $i\in\{1, \ldots, r_o\}$} ,\\
		F_{c, i, t}=&F_{c, i, t}(\varepsilon_{f_c, i, t}, \varepsilon_{f_c, i, t-1}, \ldots, ),~~\{\varepsilon_{f_c, i, t}\}_{t\in\mathbb{Z}} \text{~ is an \mbox{i.i.d.} sequence for each $i\in\{1, \ldots, r_c\}$},\\
		\me_{i, t}=&\me_{i, t}(\varepsilon_{e, i, t}, \varepsilon_{e, i, t-1}, \ldots, ), ~~\{\varepsilon_{e, i, t}\}_{t\in\mathbb{Z}} \text{~ is an \mbox{i.i.d.} sequence for each $i\in\{1, \ldots, p\}$}.
	\end{aligned} 
\end{equation}
This means that each element of vectors $\fF_{o, t}$, $\fF_{c, t}$, and $\fme_t$ is a causal process. We also assume that each process is stationary.  Let $\fmu_{\fF_o}$ be the mean of $\fF_{o, t}$.
Let $\fSigma_{\fF_o, 0}$, $\fSigma_{\fW, 0}$, and $\fSigma_{\fme, 0}$ be the covariance matrices of $\fF_{o, t}$, $\fW_t$, and $\fme_t$, respectively. We make the following assumptions.

\begin{myassump}{A}\label{ass:A} 
 Assume that (1) for any $i\in\{1, \cdots, p\}$, the processes $\{\me_{i, t}\cdot\fW_t^\top\fSigma_{\fW, 0}^{-1}\fgamma \}_{t\in\mathbb{Z}}$ and $\{\me_{i, t}\cdot(\fF_{o, t}-\fmu_{\fF_o})^\top\fSigma_{\fF_o, 0}^{-1}\fmu_{\fF_o}\}_{t\in\mathbb{Z}}$ are stationary processes. For any finite dimensional subvector $\breve{\fme}_t$ of 
$(\varepsilon_{f_o, 1, t}, \cdots, \varepsilon_{f_o, r_o, t}, \varepsilon_{f_c, 1, t}, \cdots, \varepsilon_{f_c, r_c, t}, \varepsilon_{e, 1, t}, \cdots, \varepsilon_{e, p, t})$, let $\{\breve{\fme}_t\}_{t\in\mathbb{Z}}$ be an \mbox{i.i.d.} vector sequence.  
 (2) $\mathbb{E}(\fme_t\fF_{o, t}^\top)=\mathbf{0}$ and $\mathbb{E}(\fme_t\mid \fW_t)=\mathbf{0}$; (3) $\|\fSigma_{\fF_o, 0}\|<\infty$, $\|\fSigma_{\fme, 0}\|<\infty$, and $\|\fSigma_{\fW, 0}\|<\infty$.
 \end{myassump}


 \begin{myassump}{B}\label{sparse:moment}
Suppose that $\max\{\sup_{1\leq i\leq p}\mathbb{E}(|\me_{i, t}|^d), \sup_{1\leq i\leq r_c}\mathbb{E}(|W_{i, t}|^d), \sup_{1\leq i\leq r_o}\mathbb{E}(|F_{o, i, t}|^d)\}<\infty$, for some $d>8$.
\end{myassump}

\begin{myassump}{C}\label{factorloadings}
Assume that 
there exist two constants $0<c<C<\infty$, such that $cp\leq \lambda_{\min}(\fB_c^\top\fB_c)\leq \lambda_{\max}(\fB_c^\top\fB_c)\leq Cp$.
\end{myassump}

\begin{myassump}{D}\label{ass:uncor}
Let $\mathcal{S}_i=\big\{j=1, \cdots, p: \mathbb{E}(\me_{i, t}\me_{j, t}|\fW_t)\neq\mathbb{E}(\me_{i, t}|\fW_t)\mathbb{E}(\me_{j, t}|\fW_t) \big\}$ and $N_p^*=\max_{1\leq i\leq p}\mathrm{Card}(\mathcal{S}_i)$. Assume that $N_p^*=o(p/n)$.
\end{myassump}

\begin{myassump}{E}\label{sparse:alpha}
Assume that $\| \fmu \|=o(p^{1/2}n^{-1/2}).$
\end{myassump}

\begin{myassump}{F}\label{ass:autocor_2}
Let $\{\tau_n\}$ be a positive sequence that satisfies  $\tau_n = o(n^{1/4})$ and $\tau_n/\log n\rightarrow\infty$. Assume that $pn^{1-d/4}=o(1)$ and $$\max\Big\{
p^2\cdot\sup_{1\leq i\leq p}\Theta_{\me_{i, d}}^{d/2}( \lfloor \tau_n \rfloor  ), p\cdot\sup_{1\leq i\leq r_c}\Theta_{W_{i, d}}^{d/2}( [\tau_n] ), p \cdot\sup_{1\leq i\leq r_o}\Theta_{F_{o, i, d}}^{d/2}(  [\tau_n]  )\Big\}=o(\tau_n^{d/2}).$$
\end{myassump}

\begin{rmk}
Assumption \ref{ass:A} states that idiosyncratic errors are uncorrelated with observed and latent factors, and that these three components are generated from causal processes. 
The finite-moment condition in Assumption \ref{sparse:moment} is quite mild.
Assumption \ref{factorloadings} implies that factor loadings are pervasive \citep{fan2019farmtest}.
Assumption \ref{ass:uncor} requires that the covariance matrix of idiosyncratic returns ${\fme_t}$ be sparse. As the cross-sectional correlation between the data is captured by both observed and latent factors, the sparsity assumption in Assumption \ref{ass:uncor} is reasonable; see \cite{fan2013large}. 
Assumption~\ref{sparse:alpha} implies that $\fmu$ is sparse. {This identification condition is required, unless a set of ``negative controls'' known a prior to follow the null distribution is available \citep{wang2017confounder}}.  See Section~\ref{Sec_4.2} for details.
Assumption~\ref{ass:autocor_2} implies that the dimension (p), the moments of the factors and idiosyncratic errors (d), the serial correlation, and the sample size (n) affect the convergence rate of our estimators. In general, larger dimension and stronger serial correlation result in a slower convergence rate, whereas larger moments and larger sample size lead to faster convergence.
\end{rmk}

The convergence results for $\widehat{\fB}_c$ are presented in Proposition~\eqref{lem:bhat}. Let $\ffb_{c, i}$ and $\widehat{\ffb}_{c, i}$ be the $i$-th row of $\fB_c$ and its estimator, respectively.
\begin{prop}\label{lem:bhat}
Under Assumptions \ref{ass:A}--\ref{ass:autocor_2}, 
we have 
\begin{equation}
\begin{aligned}
&\sup_{1\leq i\leq p}\Big\|\ffb_{c, i}-\frac{1}{p}\fB_c^\top\widehat{\fB}_c\widehat{\ffb}_{c, i}+\frac{1}{np}(\fH_n^{-1})^\top\widehat{\fV}_p^{-1}\widehat{\fB}_c^\top\fB_c\fW^\top \fme_{i, \cdot}^\top \Big\|
=o_p\big(n^{-1/2}\big),
\label{bhatdiff}
\end{aligned}
\end{equation}
where $\widehat{\fV}_p=\diag(\widehat{\lambda}_1^{(p)}, \ldots, \widehat{\lambda}_{r_c}^{(p)})$  
with $\widehat{\lambda}_1^{(p)}>\ldots>\widehat{\lambda}_{r_c}^{(p)}$  being the $r_c$ largest eigenvalues of 
\begin{align*}
1/(np)  \fX \qfb1 \mathcal{Q}(\widetilde{\fF}_o)\qfb1 \fX^\top,
\end{align*}
and
$\fH_n=\frac{1}{n}{\fW^\top} \qfb1\mathcal{Q}(\widetilde{\fF}_o)\qfb1\fW \times \frac{1}{p}\fB_c^\top\widehat{\fB}_c\widehat{\fV}_p^{-1}$.
\end{prop}
\begin{rmk}
From \eqref{bhatdiff}, we obtain $\big\|\ffb_{c, i}-\frac{1}{p}\fB_c^\top\widehat{\fB}_c\widehat{\ffb}_{c, i}\big\|=O_p(n^{-1/2})$, which states that $\widehat{\ffb}_{c, i}$ is a biased estimator for $\ffb_{c, i}$. This convergence result implies that $\big\|\ffb_{c, i}^\top\fgamma-\frac{1}{p}(\fB_c^\top\widehat{\fB}_c\widehat{\ffb}_{c, i})^\top\fgamma \big\|=O_p(n^{-1/2})$, from which the bias in \eqref{eq:orghat_alpha} can be eliminated. Similarly, we can show that $\|\ffb_{c, i}^\top\fW^\top (\fI_n-\widetilde{\fF}_o(\widetilde{\fF}_o^\top \widetilde{\fF}_o)^{-1}\fF_o^\top )\fb1 - \frac{1}{p}(\fB_c^\top\widehat{\fB}_c\widehat{\ffb}_{c, i})^\top \fW^\top (\fI_n- \widetilde{\fF}_o(\widetilde{\fF}_o^\top \widetilde{\fF}_o)^{-1}\fF_o^\top )\fb1 \|=O_p(n^{-1/2})$, from which the factor structure can be eliminated from \eqref{eq:orghat_alpha}.
\end{rmk}

Our main result regarding the asymptotic properties of $\widehat{\mu}_i$ is summarized in Theorem~\ref{thm:clt}.
\begin{thm}\label{thm:clt}
	Let $\widehat{\fmu}=(\widehat{\mu}_1, \ldots, \widehat{\mu}_p)^\top$. 
	If Assumptions \ref{ass:A}--\ref{ass:autocor_2} hold, as $n, p \rightarrow \infty$, we have
	\begin{align*}
		\frac{\sqrt{n}(\widehat{\mu}_i - \mu_i)}{\sigma_i}\stackrel{ \mathcal{L} }{\longrightarrow}N(0, 1)
	\end{align*}
	uniformly in $i$, where $\sigma_i^2=\sum_{t=-\infty}^\infty\mathbb{E}\Big(\me_{i, 0} \me_{i, t}FW_0  FW_t \Big)$ with $FW_t :=  1-\fW_t^\top\fSigma_{\fW, 0}^{-1}\fgamma - (\fF_{o, t}-\fmu_{\fF_o})^\top\fSigma_{\fF_o, 0}^{-1}\fmu_{\fF_o}$.
\end{thm}
\begin{rmk}
Under the static factor model, observations are \mbox{i.i.d.} sequences and $\fme_t, \fW_t$, and $\fF_{o, t}$ are mutually uncorrelated. When the factors have zero mean, i.e., $\fgamma=\mathbf{0}$ and $\fmu_{\fF_o}=\mathbf{0}$, the asymptotic variance $\sigma_i^2$ is equal to the  variance of $\fme_{i, t}$. When the means of the factors are non-zero, $\sigma_i^2$ simplifies to $\sigma_{\me, i}^2(1+ \fgamma^\top\fSigma_{\fW, 0}^{-1}\fgamma + \fmu_{\fF_o}^\top\fSigma_{\fF_o, 0}^{-1}\fmu_{\fF_o})$, which is the same as the variance in Theorem 1 of \cite{giglio2019thousands}. This implies that the means of the factors affect the standard error of the estimator. In the dynamic factor models \eqref{eq:YXZ} and \eqref{causal}, the serial correlations of $\fme_t, \fW_t$, and $\fF_{o, t}$ contribute to $\sigma_i^2$.
\end{rmk}

Now we outline the proof of Theorem \ref{thm:clt}. The details can be found in Appendix B. 
By \eqref{bernts} and Proposition \ref{lem:bhat},  we have
\begin{equation}\label{thm:alphaest}
\widehat{\mu}_i = \mu_i+\frac{\fme_{i, \cdot} \fb1 }{n}-\frac{\fme_{i, \cdot}\fW \fSigma_{\fW, 0}^{-1}\fgamma }{n}-\frac{\fme_{i, \cdot} (\fF_o-\fb1\fmu_{\fF_0}^\top) \fSigma_{\fF_o, 0}^{-1}\fmu_{\fF_o}}{n}+R_{i, n},
\end{equation}
where  $R_{i, n}$ is the remaining term. To show that $R_{i, n}=o_p(n^{-1/2})$ uniformly for $i$, we use the results 
\begin{equation*}
\sup_{1\leq i\leq p} \Big\|\ffb_{c, i}^\top\fgamma - \frac{1}{p}\big(\fB_c^\top\widehat{\fB}_c\widehat{\ffb}_{c, i}\big)^\top \fgamma + \frac{1}{np}\fme_{i, \cdot}\fW \fB_c^\top\widehat{\fB}_c\widehat{\mathbf{V}}_p^{-1}\fH_n^{-1}\fgamma \Big\|
=o_p\big(n^{-1/2}\big)
\end{equation*}
in Proposition \ref{lem:bhat}. Then we use Theorem 2.2 in \cite{Jirak16} to show that for  each $i$, 
\begin{align*}
\frac{ \sqrt{n}(\widehat{\mu}_i - \mu_i)}{\sigma_{i}}\stackrel{  \mathcal{L} }{\longrightarrow}N(0, 1).
\end{align*} 
As the convergence rates strongly depend on the serial correlation of the underlying causal processes, we further confirm  that the related bound can be uniformly controlled in $i$.

\subsection{FDR Control Theory}\label{Sec_3.2}
The FDP of the procedure in \eqref{LL} evaluated for any critical value $\xi$ is defined as
\begin{align*}
\mathrm{FDP}(\xi) := \frac{ V(\xi) }{\max\{R(\xi), 1\}}, 
\end{align*}
where $V( \xi ) :=\sum_{i\in\mathcal{I}_0}\mathbbm{1}(T_i\geq \xi)$ is the number of false discoveries and $R( \xi ) := \sum_{i=1}^p\mathbbm{1}(T_i \geq \xi )$ is the number of total rejections. The FDR  of our proposed procedure is defined as
$\mathrm{FDR}( L ) := \mathbb{E}(  \mathrm{FDP}( L ))$, where $L$ is the data-driven cutoff in \eqref{LL}.  
%

To control the FDR,  observe that
\begin{align*}
\mathrm{FDP}(L) =\frac{ V(L) }{ R(L)\vee 1 } &=\frac{ 1+\sum_{i}I(T_i \le -L)  }{ R(L)\vee 1 }\times \frac{ V(L) }{ 1+\sum_{i}I(T_i \le -L) } \cr
& \le  \beta \times \frac{ \sum_{i \in \mathcal{I}_0 }I(T_i \ge L)   }{ \sum_{i \in \mathcal{I}_0 }I(T_i \le -L) }, 
\end{align*}
from which the FDP can be controlled at the $\beta$ level if the population distribution of the test statistics is symmetric under the null hypothesis, and the empirical distribution of the negative statistics and that of the positive statistics converge to their counterparts. This is achieved through a careful analysis of an upper bound for $L$ and the following uniformly convergence results.
Specifically, we show that $0 < L\le c\log(p) $ for some $c>0$ with probability tending to one and 
\begin{equation}\label{SYM} 
\sup_{ 0 \le x  \le c\log(p) }\Bigg{|} \frac{  \sum_{i \in \mathcal{I}_0  }I(T_i \ge x) }{ \sum_{i \in \mathcal{I}_0  }I(T_i \le -x)} -1 \Bigg{|}=o_p(1).
\end{equation}
The proof can be found in Lemma~\ref{lem_mod_fdr} and Lemma~\ref{lem:testtail} in Appendix C.

The next theorem reveals that the FDR of our proposed method can be controlled by a prespecified level under some additional conditions.
\begin{myassump}{G}\label{ass:fdr}
	Assume that $\sup_{1\leq i\leq p}\sum_{\ell = 1}^\infty\ell^2\big(\sup_{t\in\mathbb{Z}}\|\me_{i, t} - \me_{i, t}^{(\ell, \prime)}\|_6 \big)<\infty$, $\sum_{\ell = 1}^\infty \ell^2\big(\sum_{j = 1}^{r_c}$
	$\sup_{t\in\mathbb{Z}} \|W_{j, t} - W_{j, t}^{(\ell, \prime)}\|_6 \big)<\infty$, and $\sum_{\ell = 1}^\infty \ell^2\big(\sum_{j = 1}^{r_o}\sup_{t\in\mathbb{Z}}\|F_{o,j, t} -F_{o,j, t}^{(\ell, \prime)}\|_6 \big) < \infty$.
\end{myassump}
Similar to Assumption~\ref{ass:autocor_2}, 
Assumption~\ref{ass:fdr} is a general measure of
serial dependence, which can be verified in many cases, see \cite{Jirak16}. 
\begin{thm}\label{thm:test}
Suppose that Assumptions \ref{ass:A}--\ref{ass:fdr} 
are satisfied. Moreover, let $\mathcal{C}_{\mu}=\{i;  \sqrt{n}|\mu_i|/\sqrt{ \log(p)  } \rightarrow \infty  \}$ and $\eta_p=\mathrm{Card}(\mathcal{C}_{\mu})$. Assume that $N_p^*/\eta_p \rightarrow 0$ and $\max_i\sigma_i^2 \cdot 2\log(p/\eta_p) \le \min_i \sigma_i^2 \cdot (1/2-\epsilon)\log(n)$ for some small $\epsilon>0$. Then, we have
\begin{align*}
\limsup_{n, p \rightarrow\infty} \mathrm{FDR}(L)  \leq \beta,
\end{align*}
where $\beta$ is a prespecified level.
\end{thm}

The derivations of our FDR control theory are significantly more complex than those in \cite{zou2020new} and \cite{Du_2021},
particularly when serial correlations exist. The main idea is to demonstrate that $T_i^{(1)}T_i^{(2)}/\sigma_i^2$ converges to the product of two independent standard normal variables uniformly for $i=1, \ldots, p$. This ensures the symmetry property of the test statistic $T_i$. The derivation involves a joint analysis of the t-statistics of two dependent splits. Specifically, in the i.i.d. setting, large deviation results for the mean can be bounded by an exponential decay rate (see Lemma S.2 in \cite{zou2020new} for details), which is crucial to verify Lemma A.1 in \cite{zou2020new} that the distribution of $T_{i}^{(1)}T_{i}^{(2)}$ is asymptotically symmetric. 
However, under a causal process, $T_i^{(1)}$ and $T_i^{(2)}$ are dependent. To prove the asymptotic symmetry of the distribution of $T_i$, the Berry--Esseen bound for one-dimensional causal processes cannot be used directly. Although \cite{gotze1983asymptotic} provide a multivariate version of the Berry--Esseen bound for weakly dependent stochastic processes, its conditional Cram$\acute{e}$r-condition seems too stringent to hold in practice. To tackle this problem, we first assume that the underlying data process, \mbox{i.e.} $\{\fW_t, \fF_{o, t}, \fme_{t}\}$, is an $m$-dependent sequence. In this case, $T_i^{(2)}$ can be decomposed into two parts:
one part is $o_p(1)$, and the other is independent of $T_i^{(1)}$, whose asymptotic distribution is a normal distribution. Thus, the asymptotic distribution of the test statistics is the product of two independent normal distributions. For the general case, we approximate the causal process by an $m$-dependent sequence with a properly selected $m$ value.

\section{Extension}\label{sec:extension}

\subsection{Improvement for Heterogeneous Cases}\label{Sec_4.1}
In our previously proposed method, we bypass the long-run variance estimation for each test and fortunately, this does not affect our FDR control because the nonnormalized test statistics still possess the symmetry properties under the null hypothesis.
When observations are homogeneous and have constant long-run variance, a modified procedure with a normalization step is equivalent to our proposed method. However, under heteroskedasticity, a test with a lower long-run variance will be missed with a high probability, which creates a fairness issue when comparing the effects of different tests. Incorporating a suitable estimate of long-run variance for each individual test will improve detection power in heterogeneous cases.
In this section, we provide a simple estimator.

Let $(\hat{e}_{i, t})_{p \times n} = \mathcal{Q}(\widehat{\fB}_c )\fX  \big( \fI_n- \widetilde{\fF}_o(\widetilde{\fF}_o^\top \widetilde{\fF}_o)^{-1}\fF_o^\top \big)$, which is the residual matrix of Step III in Section~\ref{sec:estalpha}. 
The estimator of the long-run variance $\sigma_i^2$ for each $\widehat{\mu}_i$ is defined as 
\begin{equation}\label{sn_alphai}
\hat{s}_i^2 = \frac{1}{n}\sum_{1\leq t_1, t_2\leq n} \phi\Big(\frac{t_1 - t_2}{\ell_n}\Big) \hat{e}_{i, t_1}\hat{e}_{i, t_2},
\end{equation}
where $\phi(\cdot)$ is a kernel function with $\phi(0)=1$, $\phi(x)=0$ if $|x|>1$. Furthermore, $\phi(\cdot)$ is required to be even and differentiable on $[-1, 1]$ and $\ell_n\rightarrow\infty$ and $\ell_n=o(n)$.  The theoretical properties of \eqref{sn_alphai} can be found in \cite{xiao2012covariance}. To verify that the refined procedure with long-run variance estimates can still control the FDR at the nominal level, we need to show the uniform consistency of the long-run estimators. This is beyond the scope of this manuscript, but could be addressed in future research.

\subsection{Bias Correction Using a Negative Control Set}\label{Sec_4.2}
When the magnitude of non-zero elements in $\fmu$ is large or $\fmu$ is not sparse, the finite-sample performance of the proposed test statistics in \eqref{alpha_formula} will be degraded due to the existence of non-ignorable bias. Specifically, the bias of the estimated $\fmu$ is approximated by $\mathcal{H}( \fB_c )\fmu$ when the latent factors are removed by $\mathcal{Q}(\widehat{\fB}_c)$, and it is negligible by Assumption~\ref{sparse:alpha}.
To remove this bias, we propose a bias correction procedure to estimate the mean using a negative control set. We develop the idea as follows.
\begin{enumerate}[i)]
\item 
Find a negative control set $\mathcal{S}$ that contains only null indices $i$ such that $\mu_i=0$ \citep{wang2017confounder}. 
For any random matrix $\fA$, define $\fA_{\mathcal{S} }$ as the submatrix with row indices $\mathcal{S}$.
Then, Equation \eqref{eq:f_oresid} restricted on the set $\mathcal{S}$ can be written as
\begin{equation}\label{eq:ext_resid}
\begin{aligned}
\fX_{\mathcal{S} }-(\widehat{\fB_{o, \mathcal{S}}+\fB_{c, \mathcal{S}} \fPsi} )\fF_o^\top  = 
&\fX_{\mathcal{S} } (\fI_n- \widetilde{\fF}_o(\widetilde{\fF}_o^\top \widetilde{\fF}_o)^{-1} \fF_o^\top )   \\
=& \{  \fB_{c, \mathcal{S}}( \fW+\f1_n\fgamma^\top)^\top +\fme_{\mathcal{S} } \}(\fI_n- \widetilde{\fF}_o(\widetilde{\fF}_o^\top \widetilde{\fF}_o)^{-1} \fF_o^\top ) .
\end{aligned}
\end{equation}


%
\item 
Compared with Equation \eqref{eq:f_oresid}, we find that the right-hand side of \eqref{eq:ext_resid} is irrelevant for $\fmu$. Running a Fama--MacBeth regression on the data $\{\fX_{\mathcal{S} }, \widehat{\fB}_{c, \mathcal{S} }\}$,  an unbiased estimator of the mean of the latent factors, i.e., $(\fW+\f1_n \fgamma^\top)^\top( \fI_n- \widetilde{\fF}_o(\widetilde{\fF}_o^\top \widetilde{\fF}_o)^{-1}\fF_o^\top)\f1_n/n$, which is not contaminated by $\fmu$, is given by
\begin{equation}
\begin{aligned}\label{new_est_muw}
(\widehat{\fB}_{c, \mathcal{S} }^\top\widehat{\fB}_{c, \mathcal{S} })^{-1}\widehat{\fB}_{c, \mathcal{S} }^\top\fX_{\mathcal{S} }
(\fI_n-\widetilde{\fF}_o(\widetilde{\fF}_o^\top \widetilde{\fF}_o)^{-1}\fF_o ^\top ) \f1_n/n.
\end{aligned}
\end{equation}
%
\item Plugging \eqref{new_est_muw} into \eqref{eq:f_oresid}, we can obtain the new estimator for $\fmu$ as
\begin{equation*}
\begin{aligned}
\widehat{\fmu}_{new}^\top = &n^{-1} \fX\big( \fI_n-\widetilde{\fF}_o(\widetilde{\fF}_o^\top \widetilde{\fF}_o^\top )^{-1} \fF_o^\top \big)  \fb1 
-n^{-1} \widehat{\fB}_c(\widehat{\fB}_{c, \mathcal{S} }^\top\widehat{\fB}_{c, \mathcal{S} })^{-1}\widehat{\fB}_{c, \mathcal{S} }\fX_{\mathcal{S} } (  \fI_n-  \widetilde{\fF}_o(\widetilde{\fF}_o^\top \widetilde{\fF}_o )^{-1}\fF_o^\top )\f1_n .
\end{aligned}
\end{equation*}
\end{enumerate}
Typically, some prior information is available to specify the negative control set. {For example, funds that are known a {\it priori} not to generate excess returns can constitute the negative control set.}
Otherwise, the negative control set can be determined by a threshold approach. According to \eqref{thm:alphaest} and Proposition \ref{Bernstein_ineq}, for some threshold $c_n\rightarrow0$, we have
\begin{align*}
\mathbb{P}\Big(\max_{\{i: \mu_i = 0\}}\big|\widehat{\mu}_i\big|\geq  c_n\Big) = o(1),
\end{align*}  
from which we approximate the negative control set as $\mathcal{S}=\{i=1, \ldots, p; |\widehat{\mu}_i | \leq c_n\}$.  
By Proposition \ref{Bernstein_ineq}, we choose $c_n = o(\tau_n/\sqrt{n})$. In our numerical results, we choose $c_n = O(\log n)$.

\section{Simulation}\label{sec:simu}

In our simulation studies, we examine the finite-sample performance of our proposed procedure in terms of FDR and power. We consider two scenarios: 1) the observations are independent, and 2) serial correlation exists among the observations. Both scenarios involve observable and latent factors, which collectively capture the factor structure. The FDR and power are calculated by averaging the proportions from $1,000$ replications.

\subsection{Simulation Settings}\label{sec:setup}

Consider the fund selection problem. We use the Fung--Hsieh seven-factor model, as introduced in \cite{fung2001risk}. Let 
\begin{align}\label{eq:dexf}
	\fX_t=&\fmu+\fB\fF_{t}+\fme_t,\quad t=1, \ldots, n,
\end{align}
where $\fX_t$ denotes the $p$-dimensional monthly excess returns of hedge funds, $\fF_t$ denotes the seven factors, $\fB \in \mathbb{R}^{p \times 7}$ denotes the factor loadings, and $\fme_t$ is the idiosyncratic error. To align with our model setting, we treat the first three variables as observable factors, 
which encompass a subset of systematic risks, and the last four variables as latent factors, which represent unobservable systematic risks.

Assume that $\fF_{t}\in\mathbb{R}^7$ has mean zero and covariance matrix $\fSigma_{\fF, 0}$, where $\fSigma_{\fF, 0}$ represents the sample covariance matrix of the Fung--Hsieh seven factors from 2010/01/01 to 2019/12/31. The data can be found at \href{https://people.duke.edu/~dah7/HFRFData.htm}
{\it https://people.duke.edu/~dah7/HFRFData.htm}. To generate $\fB_{sim}$, we follow these steps:
\begin{enumerate}[1.]
\item Using the same criteria as in Section \ref{sec:emp}, we clean the monthly hedge fund returns reported in the Lipper TASS dataset from January 2010 to December 2019. In total, there are $572$ hedge funds with complete data over this period. 
\item For each month, let the monthly return vectors be the dependent variable $\fX_t$, the Fung--Hsieh seven factors be the explanatory variable $\fF_{t}$, and assume that they satisfy Model \eqref{eq:dexf}.
\item By applying the OLS method to the data $\{\fX_t, \fF_{t}\}$, we can obtain $\widehat{\fB}$ for $\fB$ and the estimated residual vectors. \label{sim:res}
\item Let $\fmu_{B}$ and $\fSigma_{B}$ be the mean and covariance matrix of the rows of $\widehat{\fB}$, respectively. For each row of $\fB_{sim}$, we generate the data independently from $\mathcal{N}(\fmu_{B}, \fSigma_{B})$.
\end{enumerate}
The covariance matrix of the idiosyncratic error $\fme_t$ is $\fSigma_{\fme, 0}$, where $\fSigma_{\fme, 0} = \big(\rho^{|i-j|}\big)$ with $\rho=0.5$. 
For notational convenience, we refer to our method as YD. To demonstrate the merit of our proposed method, the following approaches are compared:
\begin{enumerate}[a)]
\item The factor-adjusted multiple testing procedure (FAT) proposed by \cite{lan2019factor}.
\item The factor-adjusted robust multiple testing method (FARM) \citep{fan2019farmtest}, where we use the \texttt{R} package \texttt{FarmTest}.
\item The modified mean screening BH procedure (SBH) \citep{giglio2019thousands}, where the $p$ values are obtained by normal calibration. We do not implement the screening step because we focus on two-sided alternative hypotheses.
\item The self-normalization (SN) approach proposed by \cite{shao2010self}. Specifically, after controlling for the observed and latent factors, we apply the SN approach to the residual matrix to estimate the long-run variance of the estimated $\mu_i$. The $p$ values are calibrated using a mixture of standard normal distributions.
\end{enumerate}
Let $\pi$ be the proportion of false hypotheses. The true mean is $\fmu=(\mu_1,  \ldots,  \mu_p)^\top$, with $\mu_i=\nu$, $1\leq i\leq [ \pi p/2 ]$, $\mu_i=-\nu$, $[ \pi p/2 ]+1\leq i\leq [\pi p ]$, and $\mu_i=0$ otherwise. 
For simplicity, we set $(n, p, \pi)=(200, 1000, 0.1)$, $\nu=0.2, 0.3$, and $\beta=5\%, 10\%, 15\%$, while the results for the other parameter settings are qualitatively similar.  
%
\subsection{The Case Where the Observations Are Independent}\label{sim:onlylatentiid} 
We begin by considering the case where the observations are independent over time.
To investigate the robustness of different methods, we consider two distributions for the factors and idiosyncratic errors: (i) multivariate normal and (ii) log-normal.  The covariance matrix structures for the factors and idiosyncratic errors are specified in Section~\ref{sec:setup}. 
The results for the FDR and power are presented in Table \ref{both_iid}. The following observations are made:
\begin{enumerate}[(a)]
\item
In a normal distribution, FAT and FARM control the FDR reasonably well and their power is comparable. In contrast, our method achieves better FDR control, while the power loss is relatively small due to the simultaneous use of two test statistics. To account for non-tradable factors, SBH uses a Fama--Macbeth regression based on combined factor loadings, which inflates the standard error of the estimated risk premium. As a result, its FDR is slightly out of control. Our simulation results on SN confirm that it diminishes size distortion and thus loses some power by incorporating an inconsistent estimator of long-run variance. 

\item 
In the case of log-normal distribution, FAT, SBH, and SN fail to control the FDR at the prespecified level due to their sensitivity to heavy-tailed distributions. Our approach estimates the number of false rejections using a nonparametric method, thus the FDR can be well controlled at the significance level. FARM performs best in terms of FDR and power when signal strength is large, but it is computationally intensive, compared with our method.
\end{enumerate}

\begin{table}[h]
	\begin{center}
		\ra{0.63}
		\setlength{\tabcolsep}{1.4pt}
		\begin{tabular}{@{}c|cccccccccc@{}}\toprule[1pt]
			&~~~&~~~&\multicolumn{3}{c}{FDR$(\%)$} & \phantom{~}& \phantom{~}& \multicolumn{3}{c}{power($\%$)}\\
			\cmidrule{4-6} \cmidrule{9-11}
			&~~$~~\nu~~$&~~&$~~\beta=5$~~&~~$10$~~~~~&~~~~~$15$~~& \phantom{~}& \phantom{~}  &~~$\beta=5$~~&~~$10$~~~&~~~$15$~~ \\ 
			\midrule[1pt]
			&\multicolumn{10}{c}{Both observed and confounding factors exist} \\
			\midrule[1pt]
			\multirow{8}{*}{Normal}&\multirow{4}{*}{$0.2$}&FAT&$5.13(3.19)$&$10.40(4.23)$&$15.62(4.86)$&\phantom{~}& \phantom{~}&$52.28(10.33)$&$64.77(9.30)$&$72.09(8.46)$\\
			&&SBH&$7.05(3.73)$&$13.29(4.47)$&$19.09(5.16)$&\phantom{~}& \phantom{~}&$56.72(9.54)$&$68.19(8.49)$&$74.63(7.75)$\\
			&&FARM&$4.45(3.30)$&$8.80(4.03)$&$13.11(4.62)$&\phantom{~}& \phantom{~}&$47.96(9.63)$&$59.95(8.77)$&$67.52(8.30)$\\
			&&SN&$7.28(8.94)$&$12.61(7.02)$&$17.97(7.19)$&\phantom{~}& \phantom{~}&$15.14(7.52)$&$26.23(9.41)$&$35.62(10.50)$\\
			&&YD&$3.81(4.10)$&$8.70(5.35)$&$13.45(6.19)$&\phantom{~}& \phantom{~} &$39.29(19.97)$&$56.32(13.79)$&$64.87(11.74)$\\
			\cmidrule{2-11}
			&\multirow{4}{*}{$0.3$}&FAT&$5.34(2.46)$&$10.86(3.56)$&$16.29(4.46)$&\phantom{~}& \phantom{~}&$96.04(2.61)$&$97.86(1.78)$&$98.61(1.38)$\\
			&&SBH&$7.42(2.95)$&$13.73(4.01)$&$19.97(4.73)$&\phantom{~}& \phantom{~}&$96.67(2.28)$&$98.12(1.69)$&$98.77(1.31)$\\
			&&FARM&$4.46(2.47)$&$8.83(3.52)$&$13.18(4.42)$&\phantom{~}& \phantom{~}&$94.87(2.88)$&$97.15(2.00)$&$98.22(1.49)$\\
			&&SN&$6.26(3.38)$&$12.03(4.17)$&$17.78(4.66)$&\phantom{~}& \phantom{~}&$58.67(8.17)$&$71.87(7.30)$&$79.75(6.27)$\\
                 &&YD&$4.18(3.03)$&$9.05(4.46)$&$14.14(5.42)$&\phantom{~}& \phantom{~}  &$92.95(4.58)$&$96.27(2.78)$&$97.51(2.09)$\\
         \midrule[1pt]
         \multirow{8}{*}{log-normal}&\multirow{4}{*}{$0.2$}&FAT&$13.20(4.65)$&$19.38(4.76)$&$24.52(5.17)$&\phantom{~}& \phantom{~}&$59.32(10.19)$&$70.97(8.69)$&$77.30(7.55)$\\
         &&SBH&$14.78(5.17)$&$20.85(5.63)$&$26.28(5.83)$&\phantom{~}& \phantom{~}&$64.86(8.71)$&$74.56(7.52)$&$80.10(6.33)$\\
			&&FARM&$8.51(8.83)$&$13.02(10.25)$&$16.95(10.97)$&\phantom{~}& \phantom{~}&$43.27(12.60)$&$54.46(12.71)$&$61.62(12.21)$\\
			&&SN&$12.92(7.91)$&$19.22(7.16)$&$24.77(7.12)$&\phantom{~}& \phantom{~}&$20.25(7.83)$&$32.59(9.14)$&$42.23(9.89)$\\
			&&YD&$4.41(4.61)$&$9.00(5.84)$&$13.60(6.54)$&\phantom{~}& \phantom{~} &$34.39(21.81)$&$53.87(16.15)$&$64.13(12.32)$\\
			\cmidrule{2-11}
			&\multirow{4}{*}{$0.3$}&FAT&$11.30(3.50)$&$17.74(4.21)$&$23.19(4.68)$&\phantom{~}& \phantom{~}&$94.29(2.87)$&$95.96(2.44)$&$96.84(2.12)$\\
			&&SBH&$12.31(4.08)$&$18.75(4.70)$&$24.75(5.26)$&\phantom{~}& \phantom{~}&$94.61(2.78)$&$96.05(2.37)$&$96.93(2.12)$\\
			&&FARM&$7.41(8.25)$&$11.94(9.74)$&$16.22(10.72)$&\phantom{~}& \phantom{~}&$89.48(7.31)$&$93.07(5.41)$&$94.89(4.43)$\\
			&&SN&$9.82(3.89)$&$16.22(4.68)$&$22.41(5.32)$&\phantom{~}& \phantom{~}&$62.65(7.78)$&$74.92(7.09)$&$82.22(6.22)$\\
         &&YD&$4.60(3.30)$&$9.23(4.64)$&$14.24(5.45)$&\phantom{~}& \phantom{~}  &$92.10(5.69)$&$95.74(2.99)$&$97.03(2.30)$\\
			\bottomrule[1pt]
		\end{tabular}
	\end{center}
	\caption{Empirical FDR and power of FAT, SBH, FARM, SN, and our method (YD) at prespecified FDR levels $\beta=0.05, 0.1$, and $0.15$, for varying signal strengths $\nu$. The data are generated from a normal distribution or a log-normal distribution. }\label{both_iid}
\end{table}

\subsection{The Case Where the Observations Are Time-dependent}\label{sim:timedep}
We next consider the case when serial correlation exists between factors and idiosyncratic errors. 
To simulate the factors, we first fit the Fung--Hsieh seven factors using GARCH(1, 1) models, from which we generate the simulated factors. 
The first three variables serve as observable factors.

To simulate the idiosyncratic errors, we analyze the residuals obtained by fitting Model \eqref{eq:dexf} with monthly hedge fund returns as observations and the Fung--Hsieh seven factors as factors; the dataset is used in Section~\ref{sec:emp}. Using the Breusch--Godfrey test, we find that the monthly returns of $33.4\%$ of  hedge funds are serially correlated. For each of these hedge funds, we obtain its residuals and fit the ARMA($k_1, k_2$) models, and collectively the residuals are generated from a mixture of ARMA models with eight components. The mixing proportion for each ARMA($k_1$, $k_2$) model is its relative frequency, while the corresponding parameters are simply set as the average of the parameter estimates for each component. We generate simulated idiosyncratic errors from the mixture distribution. Overall, the percentage of simulated idiosyncratic errors based on ARMA models is $33.1\%$. The remaining $66.9\%$ of simulated errors are sampled from the normal distribution. For each time point $t$, random errors are standardized and multiplied by the square root of the covariance matrix $\fSigma_{\fme, 0}$.

Table~\ref{both_ts} presents the FDR and power in the presence of serial correlations. The FAT, SBH, and FARM methods are only valid when independent, which explains why their FDRs are largely out of control. The SN method is designed for serially correlated observations and should deliver satisfactory performance on FDR control, but its power is considerably sacrificed. In comparison, our method demonstrates excellent performance on FDR control and maintains reasonably high power.
\begin{table}[h]
	\begin{center}
		\ra{0.63}
		\setlength{\tabcolsep}{1.4pt}
		\begin{tabular}{@{}cccccccccc@{}}\toprule[1pt]
			~~~&~~~&\multicolumn{3}{c}{FDR$(\%)$} & \phantom{~}& \phantom{~}& \multicolumn{3}{c}{power($\%$)}\\
			\cmidrule{3-5} \cmidrule{8-10}
			~~$~~\nu~~$&~~&$~~\beta=5$~~&~~$10$~~~~~&~~~~~$15$~~& \phantom{~}& \phantom{~}  &~~$\beta=5$~~&~~$10$~~~&~~~$15$~~ \\ 
			\midrule[1pt]
			\multicolumn{10}{c}{Both observed and confounding factors exist} \\
			\midrule[1pt]
			\multirow{5}{*}{$0.2$}&FAT&$11.09(4.60)$&$18.29(5.42)$&$24.86(5.76)$&\phantom{~}& \phantom{~}&$54.59(9.39)$&$66.76(8.49)$&$73.95(7.63)$\\
			&SBH&$13.80(4.70)$&$21.80(5.36)$&$28.99(5.71)$&\phantom{~}& \phantom{~}&$59.20(8.68)$&$69.98(7.91)$&$76.32(7.21)$\\
			&FARM&$8.51(8.83)$&$13.02(10.25)$&$16.95(10.97)$&\phantom{~}& \phantom{~}&$43.27(12.60)$&$54.46(12.71)$&$61.62(12.21)$\\
			&SN&$7.36(8.42)$&$12.82(7.63)$&$18.27(7.68)$&\phantom{~}& \phantom{~}&$13.05(6.90)$&$22.96(8.89)$&$31.88(9.95)$\\
                        &YD&$4.25(5.04)$&$9.50(6.19)$&$14.50(6.55)$&\phantom{~}& \phantom{~}&$26.34(20.30)$&$46.26(15.70)$&$57.04(12.10)$\\		
			\cmidrule{2-10}
			\multirow{4}{*}{$0.3$}&FAT&$10.19(3.57)$&$17.72(4.53)$&$24.64(5.16)$&\phantom{~}& \phantom{~}&$95.45(2.81)$&$97.51(1.96)$&$98.42(1.49)$\\
			&SBH&$12.87(3.84)$&$21.32(4.67)$&$28.96(5.11)$&\phantom{~}& \phantom{~}&$96.08(2.61)$&$97.90(1.74)$&$98.62(1.34)$\\
			&FARM&$7.41(8.25)$&$11.94(9.74)$&$16.22(10.72)$&\phantom{~}& \phantom{~}&$89.48(7.31)$&$93.07(5.41)$&$94.89(4.43)$\\
                         &SN& $6.48(3.48)$&$12.19(4.37)$&$17.90(4.98)$&\phantom{~}& \phantom{~}&$54.36(8.03)$&$67.84(7.41)$&$76.24(6.50)$\\
                      &YD&$4.71(3.42)$&$9.57(4.55)$&$14.70(5.93)$&\phantom{~}& \phantom{~} &$88.72(6.08)$&$93.68(3.86)$&$95.61(2.98)$\\
        	\bottomrule[1pt]
		\end{tabular}
	\end{center}
	\caption{Empirical FDR and power of FAT, SBH, FARM, SN, and our method (YD) at prespecified FDR levels $\beta=0.05, 0.1$, and $0.15$, for varying signal strengths $\nu$ when serial correlation exists (i.e., factors are simulated from GARCH(1, 1) and idiosyncratic errors are generated from a mixture of ARMA models.).  }\label{both_ts}
\end{table}

\subsection{The Heterogeneous Case}
\begin{figure}[h]
	\begin{center}
	\includegraphics[width=0.8\textwidth]{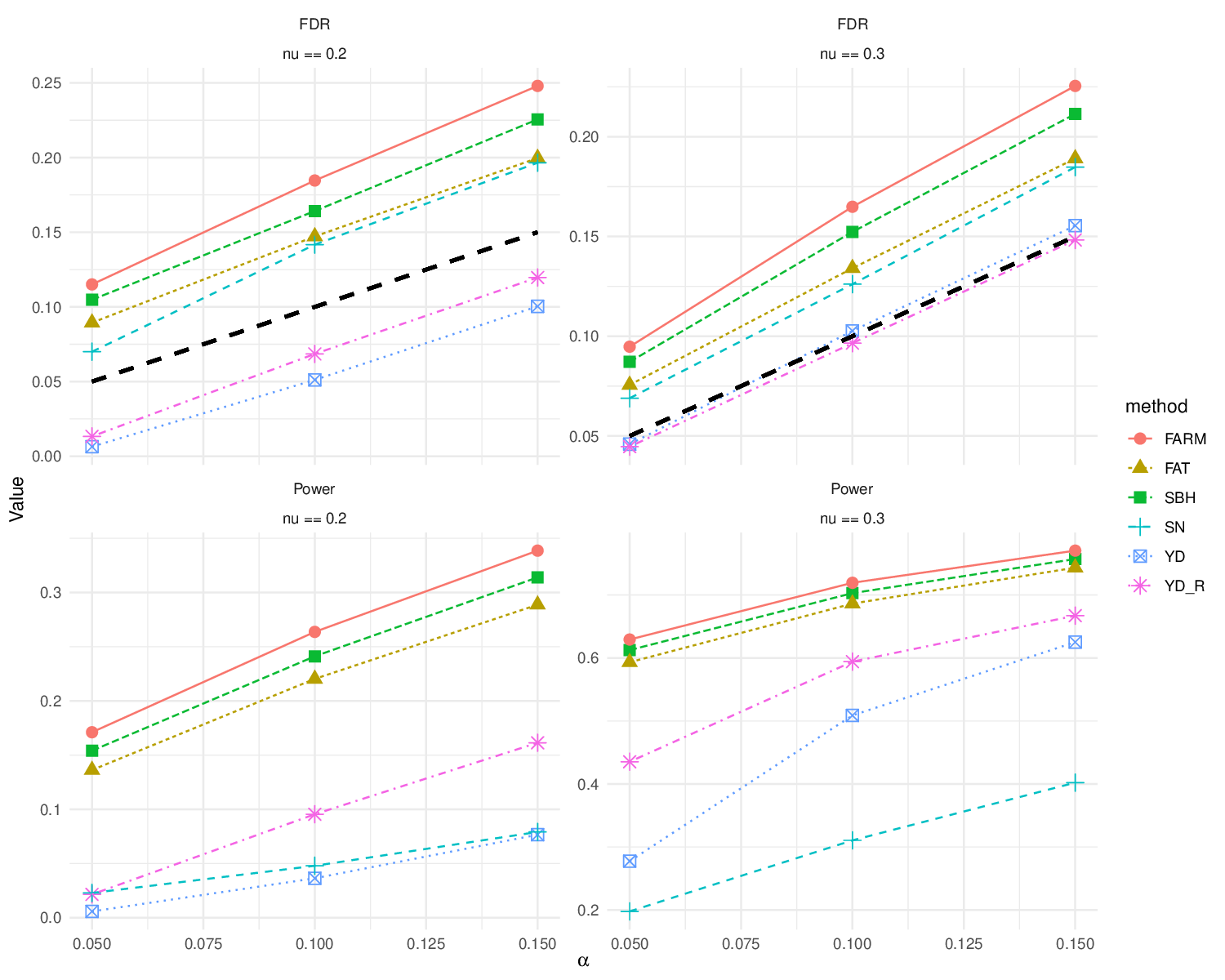}
	\caption{Empirical FDR and power of our method (YD), its improvement (YD$\_$R), and competitors (i.e., FARM, FAT, SBH, SN) at prespecified FDR levels $\beta=0.05, 0.1$, and $0.15$, for varying signal strengths $\nu$ when serial correlation exists and idiosyncratic errors are heterogeneous. }\label{IMPlatent_ts}
	\end{center}
\end{figure}

For this case, the data are generated using the same procedure as in Section \ref{sim:timedep}, except that the diagonal elements of $\fSigma_{\fme, 0}$ are sampled independently  from the uniform distribution $[1, 3]$. We compare our proposed procedure with the refined method introduced in Section \ref{Sec_4.1}, referred to as YD$\_$R, along with its competitors (i.e., FARM, FAT, SBH, SN). To estimate the long-run variance of the refined method, we adopt the Bartlett kernel, i.e., $\phi(x) = (1- |x|)\mathbbm{1}( |x|\leq 1 )$, and the bandwidth $\ell_n$ is set to $n^{1/5}$. To ensure a fair comparison, all competitors are equipped with the Heteroskedasticity and Autocorrelation Consistent (HAC) estimator for long-run variance. We apply this approach to the residuals using the \texttt{vcovHAC} function from the R package \texttt{sandwich}.
The comparisons of FDR and power for YD, YD$\_$R, and competitors are presented in Figure~\ref{IMPlatent_ts}. For different signal magnitudes and FDR levels, YD$\_$R takes into account the long-run variance estimate, so its power is boosted dramatically compared with our initially proposed method. With the data-driven choice of long-run variance, the FDR of YD$\_$R is still very close to the nominal level. In comparison, all of the competitors exhibit inflated FDRs, even when incorporating long-run variance.

\section{Empirical Analysis}\label{sec:emp}


Distinguishing skilled fund managers from others not only helps investors achieve consistently superior performance but is also a key topic in asset management. As the abnormal expected return relative to a benchmark model reflects fund managers' ability to exploit arbitrage opportunities, the multiple testing problem in \eqref{hypo} can be used to identify funds that exhibit positive or negative risk-adjusted returns (i.e., alpha).  In other words, skilled funds have managers with sufficient stock-picking capacity to cover trading costs and expenses (i.e., $\mu_i>0$), while unskilled funds do not (i.e., $\mu_i<0$). The remaining funds are classified as zero-alpha funds (i.e., $\mu_i=0$). Based on the Berk and Green equilibrium \citep{berk2004mutual}, only a very small proportion of $\alpha_i$s are expected to be non-zero.
Seminal works include \cite{barras2010false}, \cite{fama2010luck} and \cite{harvey2016and}. 

This section investigates the performance of thousands of hedge fund managers. We use monthly hedge fund returns reported in the Lipper TASS dataset from January 2005 to December 2019. To control for backfill bias,  it is recommended to remove returns prior to the date of addition to the dataset.  
We apply the technique proposed by \cite{jorion2019fix} to infer the dates of addition when they are not available. Funds that do not publish monthly returns or provide net-of-fee returns are excluded. We also remove funds whose returns are less than $-100\%$ or greater than $200\%$, and drop funds when their AUM is less than  $\$$10m, which are more likely to manipulate reporting on the dataset. 
Furthermore, we require funds to report returns for more than $5$ years without missing data. There are $2,550$ hedge funds in our dataset.

 We first verify the existence of serial correlation in their returns. After regressing the monthly returns of each hedge fund on the Fung--Hsieh seven factors, we test for the presence of serial correlation in the residuals using the  Breusch–Godfrey test. If we set the significance level at $10\%$, the empirical size is $33.02\%$. The severely distorted empirical size implies that hedge fund returns are serially correlated. Such serial correlation can render methods based on independent observations ineffective; in contrast, it is perfectly compatible with our  methodology.  
 
We start with portfolio construction based on our method, as well as the PCA\footnote{Because there are less than seven significant factors based on Eigendecomposition, we only consider the case of latent factors.} method and the individual-test method (IND). 
For each method, we select a fund at the beginning of each month, when two conditions are both satisfied: 1) its estimated alpha is positive; 2) it is significant according the corresponding test. The estimated alpha and the test result are obtained based on the monthly returns of the previous five years. Next, we invest in the selected funds in equal proportions. For our method and the PCA method, the FDRs are set to $5\%$, $10\%$, and $15\%$. Correspondingly, we reject the null hypothesis by setting the significance levels to $5\%$, $10\%$, and $15\%$ for the IND. We refer to these portfolios as the construction method.

Table \ref{tab:sharpe} reports the Sharpe ratio and the average monthly returns of the IND, PCA, YD, and YD$\_$TH portfolios from January $2010$ to December $2019$. YD$\_$TH represents the threshold method proposed in Section \ref{Sec_4.2}

\begin{table}[h]
	\begin{center}
		\ra{0.63}
		\setlength{\tabcolsep}{1.4pt}
		\begin{tabular}{@{}ccccccccc@{}}\toprule[1pt]
			~~~&\multicolumn{3}{c}{Sharpe ratio} & \phantom{~}& \phantom{~}& \multicolumn{3}{c}{Average monthly returns}\\
			\cmidrule{2-4} \cmidrule{7-9}
			~~& ~~~~FDR$=5\%$~~~~& $10\%$~~~~&  $15\%$ ~~~~& \phantom{~}& \phantom{~}  &~~~~FDR=$5\%$~~~~&$10\%$~~~~&$15\%$ ~~~~\\ 
			\midrule[1pt]
			IND& $0.88$ & $0.82$~~~~& $0.78$~~~~& \phantom{~}& \phantom{~}  & $0.61$ & $0.61$~~~~& $0.60$~~~~\\
			PCA& $0.91$ & $0.81$~~~~& $0.75$ ~~~~& \phantom{~}& \phantom{~}  & $0.61$ & $0.59$~~~~& $0.58$~~~~\\
			YD& $1.03$ & $1.13$~~~~& $1.01$~~~~&\phantom{~} & \phantom{~}  & $0.53$  & $0.58$ ~~~~& $0.56$~~~~\\
			YD$\_$TH& $1.06$& $1.14$~~~~& $1.01$~~~~& \phantom{~}& \phantom{~}  & $0.55$ & $0.58$~~~~& $0.56$~~~~\\
			\bottomrule[1pt]
		\end{tabular}
	\end{center}
	\caption{The Sharpe ratio (left) and average monthly returns (right) of the IND, PCA, YD, and YD$\_$TH portfolios from January $2010$ to December $2019$.}\label{tab:sharpe}
\end{table}

From Table \ref{tab:sharpe}, we make the following observations: 
\begin{enumerate}[1)]
\item Although the average monthly returns of the IND portfolio are slightly higher than those of the PCA and YD portfolios,  the Sharpe ratio of the IND portfolio is lower than that of the other methods. The reason may be that if we set the significance level at the same level as the FDR for each individual test, there is a high probability that the true FDR will be larger than in other methods. Thus, more underperforming hedge funds will be included in the portfolio, leading to greater volatility without any noticeable improvement in returns. 
\item The average monthly returns of the PCA and YD portfolios are comparable. 
The Sharpe ratio of the PCA portfolio is higher than that of the IND portfolio, but lower than that of the YD portfolio. As the monthly returns of around $33\%$ of hedge funds are serially correlated, the null distribution of the test statistics used in the PCA method is no longer correct. This could make empirical FDRs higher than the nominal level, and thus lead to more volatile portfolios. 
\item
The performance of YD$\_$TH and YD is comparable, implying that hedge fund signals are weak and sparse.
\end{enumerate}

\section{Conclusions}\label{sec:conclude}
We conclude this manuscript by identifying two directions for future research.

First, we incorporate long-run variance estimates in Section~\ref{Sec_4.2} and numerically demonstrate that such a refined method can substantially improve power while maintaining the FDR. It would be interesting to provide some theoretical guarantees on such long-run variance estimates so that the FDR can be controlled under heteroskedasticity. 

Second, it would be practically useful to adapt our method to the streaming process, where data are updated sequentially and decisions need to be made once new data arrive. A new sample-splitting scheme is required because signals from sequential observations change dynamically over time.


\begin{spacing}{1}
\bibliographystyle{asa}
\bibliography{reference}
\end{spacing}

\end{document}